%\divide\hsize by 1414
%\multiply \hsize by 1000
%\divide\vsize by 1414
%\multiply \vsize by 1000
\overfullrule=0pt

\def\frac#1#2{{#1\over #2}}
\def\textit#1{{\it #1\/}}
\def\CPI{{}\ifmmode \mathop{\rm CPI}\nolimits\else $\mathop{\rm CPI}\nolimits$ \fi}
\def\CPIR{{}\ifmmode \mathop{\rm CPIr}\nolimits\else $\mathop{\rm CPIr}\nolimits$ \fi}
\def\D{{}\ifmmode \mathop{\rm D}\nolimits\else $\mathop{\rm D}$ \fi}

\def\Cituj#1{\def\Znacka{\csname KlicoveSlovo#1\endcsname}%
\if\Znacka\relax ??.??.??\else \Znacka\fi}

\def\Oddelovac{
\medskip
\centerline{\hbox to .3333\hsize{\hrulefill}}
\bigskip
}

\newcount\NofEq
\def\Nq{\global\advance\NofEq by 1\hfill\hbox %to 0 pt
{(\number\NofEq)\hss}}

\def\NEq{\global\advance\NofEq by 1\relax\eqno{(\number\NofEq)}}

\def\Cr{\Nq\cr}

\def\HodnotaAktualnihoCisla{\number\NofEq}

%-cros-ref-% \immediate\openin10 \jobname.citace
%-cros-ref-% \ifeof10 \immediate\openout10\jobname.citace\immediate\closeout10\fi
%-cros-ref-% \immediate\closeout10
%-cros-ref-% \immediate\closein10
%-cros-ref-% 
%-cros-ref-% \input \jobname.citace
%-cros-ref-% \immediate\openout10=\jobname.citace
%-cros-ref-% 

\def\PripravCitaci#1{\global\edef\HodnotaAktualnihoCisla{\number\NofEq}%
%-cros-ref-% \immediate\write10{\noexpand\expandafter\noexpand\gdef\noexpand\csname
%-cros-ref-% KlicoveSlovo#1\noexpand\endcsame{\HodnotaAktualnihoCisla}}%
\expandafter\xdef\csname KlicoveSlovo#1\endcsname{\HodnotaAktualnihoCisla}}

\long\def\ABSTRAKT#1\Oddelovac{{
\leftskip 2\parindent
\rightskip2\parindent
\parindent=0pt
\bigskip
\noindent{\bf Abstract:} #1
\bigskip
}}

\def\KEYWORDS{\noindent{\bf Keywords:} }
\def\MSC{\noindent{\bf MSC 2000 Classification:} }

\newcount\nrA
\newcount\nrB

\def\K#1{\medskip\advance\nrA by 1 \noindent{\bf\number\nrA. 
 #1:}}

\def\EXAMPLE{%\K{Example #1}
\medskip\advance\nrB by 1 \noindent{\bf\number\nrA.\number\nrB. Example}
}

\def\PROOF{\advance\nrB by 1\smallskip\noindent{\bf\number\nrA.\number\nrB. Proof:\/}}

\def\DEF{\advance\nrB by 1\smallskip\noindent{\bf\number\nrA.\number\nrB. 
Definition: }}

\def\THEOREM{\advance\nrB by 1\smallskip\noindent{\bf\number\nrA.\number\nrB. 
Theorem: }}

\def\LEMMA{\advance\nrB by 1\smallskip\noindent{\bf\number\nrA.\number\nrB. 
Lemma: }}

\def\REFERENCES{%
\bigskip
\parindent=0pt
{\bf References:}} 

\centerline{\bf Functional Equation of the Rate of Inflation}
\smallskip 
\centerline{V clav Studenì}
\bigskip
\centerline{%
Department  of Applied  Mathematics, Faculty of Economics}  
\centerline{ and Administration, Masaryk University}
\centerline{Lipov  24a, 602 00 Brno}
{\catcode`\@=11
\centerline{studeny@econ.muni.cz}
}
%\footnote{}{The recherche is supported  by\dots}
\ABSTRAKT
This short note aims to  introduce  a rule which admits to compute %any time rate of 
interest in any time per any time, 
%or any time 
rate of inflation per any time in any moment, if  the rate 
of interest or the rate of inflation by unity of time is an arbitrary  integrable 
function.
 
The main result  is the generalization of the well known rule 
$\iota(t)+1=\frac {{\rm X}(0)}{{\rm X}(t)}=\prod _{i=0}^{n - 1}\,({I_{i}} + 
1)^{({t_{i + 1}} - {t_{i}})}$%
, 
which holds for the piecewise constant approximation of the rate of inflation with 
the values $I_i$ to rule
{%
$\iota(t)+1=\frac {{\rm X}(0)}{{\rm X}(t)}=e^{ \left(  \! \int _{0}^{t}{\rm 
ln}(1 + \iota (u))\,\mathop{\rm d}u
 \!  \right) }$%
}, 
which should be used for arbitrary integrable function $\iota$.

The usage of the rule is demonstrated on the examples based on real data of index of 
prices of non regulated prices in Czech republic.

\KEYWORDS interest rate, rate of inflation, functional equations, integral equation, 
approximation 

\MSC
62P05 Applications to actuarial sciences and financial mathematics,
% 65C30 Stochastic differential and integral equations
65R20 Integral equations,
39B72 Systems of functional equations and inequalities.

\Oddelovac

\K {Introduction}
Let $\CPI_i(t)$ be the  index of prices in time $t$, let $X(t)$ be the value of an 
objective function (the real value of  a unit  of money), 
let $\iota(t)$ be rate of inflation per time $\langle0,1\rangle$. Then 
$$
{X(t_0)\over X(t_1)}=
{\CPI(t_1)\over 
\CPI(t_0)}=1+\iota(\langle t_0,t_1\rangle)\NEq$$

We can measure  \CPI in any time,  consequently  we  can fix $\iota(\langle 
t_0,t_1\rangle)$ in an arbitrary interval $\langle t_0,t_1\rangle$

\DEF If $\iota(\langle t_0,t_1\rangle)$ in some interval 
$\langle T_0,T_1\rangle\supset\langle t_0,t_1\rangle$ 
depends only on the length  $t_1-t_0$ of the interval $\langle t_0,t_1\rangle$ but not on 
the origin $t_0$ of the interval we call the
inflation constant (on interval $\langle T_0,T_1\rangle$).

\LEMMA The inflation is constant, 
if and only if  the \CPI is exponential function
$x\mapsto a\cdot e^{bt}$.

Let the inflation is constant and let us denote
$\iota(t_1-t_0)=\iota(\langle t_0,t_1\rangle)$, then it holds
\def\Cr{\Nq\cr}
$$\displaylines{
\hfill1+\iota(t_i)=(1+\iota(t_j))^{t_j\over t_i}\hfill\Cr
\noalign{\leftline{especially}}
\hfill 1+\iota(t)=(1+\iota(1))^{t}\hfill\Cr}
$$
and the real value (value of the objective function) in the time $t$ is
$$
X(t)={X(0)\over (1+\iota)^t}
\NEq$$

Let the inflation depends differentiably on the time in a time $t_0$ almost 
everywhere. We define the rate of
inflation per unity of time in time $t_0$ as rate of constant inflation which has first 
order contact with the inflation in question almost everywhere and on the rest of the 
points we can define it by one sided limits. 

More precisely:
The solution of equations:
$$\eqalign{
F(x)=Ae^{Bx}&=\CPI(x)\cr
\D(F)(x)=\D\left(A{e^{Bx}}\right)=AB{e^{Bx}}&=\D(\CPI)(x)
}
\NEq$$
is
$$
\eqalign{
A&=\CPI(x)\left ({e^{{\frac {\left (\D (\CPI)(x)\right )x}{
\CPI(x)}}}}\right )^{-1}
\cr
B&={\frac {{\partial\over\partial x} (\CPI)(x)}{\CPI(x)}}
}
\NEq$$
Consequently
$$\iota=\iota(1)=
x\mapsto {F(x+1)\over F(x)}-1
=x\mapsto {e^{{\frac {\D (\CPI)(x)}{\CPI(x)}}}}-1
=x\mapsto {e^{\D (\ln\circ\CPI)(x)}}-1   
\NEq$$
\DEF
Let \CPI be the index of prices, then the mapping
$$
x\mapsto {e^{{\frac {\D (\CPI)(x)}{\CPI(x)}}}}-1
\NEq
$$
is the rate of inflation per unity of time (in time $x$).

\K {Functional equation of rate of inflation}

%%%

% We can substitute in any point $t$ by constant inflation with the same 
% rate per wery small period of time more precisely:

The most general case considered so far is that one  with
 $\CPI_i$ piecewise  exponential function $x\mapsto a_ie^{b_ix}$ with the constant 
$a_i$, $b_i$ on an interval $\langle t_i,t_{i+1}\rangle$, $(t_0=0)$, then 
the rate of inflation per unity of time is
piecewise  constant  function and the rate of 
inflation per time $\langle 0,t\rangle$ is 
$\prod _{\{i|t_i\leq t\}} {\CPI_{i+1}\over\CPI_i}-1$

Let us suppose that $\CPI$ is the integrable function, 
$\iota$ is the rate of inflation per unity of time (integrable function), $X$ the 
objective function (the opposite value of $\CPI$ in the suitable rescaling)

If an inflation is constant and if its rate of time of duration $t$  is equal to
$\iota (t)$, then
$$1 + \iota (t)
=
(1+\iota (1))^{t}
=
\frac {{\rm X}(0)}{{\rm X}(t)}.
\NEq$$%

If an inflation is piecewise constant (the rate of the  inflation per  time unit is 
piecewise constant function) and if the rate of  inflation per a unit of time in 
every point of interval $%{I_{i}}=
\langle {t_{i}}, \,{t_{i + 1}}\rangle \,;\,i=0 \dots n$
is equal
${I_{i}}$%
,then for the rate of inflation of time $t$
{%
$\langle {t_0}, \,{t_n}\rangle $%
} holds   
{%
$$%1 + \iota (t)
%=
1+\iota (\langle t_0, \,t_n\rangle )=
(\prod _{i=0}^{n - 1}\,({I_{
i}} + 1)^{({t_{i + 1}} - {t_{i}})})
=
\frac {{\rm X}(t_0)}{{\rm X}(t_n)}.
\NEq$$%
}

This note aims to find the rate of inflation in the time $t$ per a time period
{%
$\langle {t_{0}}, \,{t_{1}}\rangle $%
} (
{%
$\langle 0, \,t\rangle $%
}, after   suitable  rescaling), 
in the case that the rate of inflation per unit of time  is an arbitrary integrable 
function.

\THEOREM
If the rate of   inflation per unit of time 
$t$
equals
{%
$\iota (t)$%
}, 
then the rate of inflation per time
$\langle 0, \,t\rangle $
 (plus one) is equal to
{%
$$\iota (\langle t_0, \,t_n\rangle )
=
\frac {{\rm X}(t_0)}{{\rm X}(t_n)}-1
=
e^{ \left(  \! \int _{t_0}^{t_n}{\rm 
ln}(1 + \iota (u))\,\mathop{\rm d}u
 \!  \right) }
-1
\NEq$$%
}
and the equation 
{%
$$\frac {{\rm X}(t_0)}{{\rm X}(t_n)}=\prod _{i=0}^{n
 - 1}\,({I_{i}} + 1)^{({t_{i + 1}} - {t_{i}})}\NEq$$%
}
is a  special case of the previous one  for the piecewise constant function $\iota$ with 
values $I_i$ on intervals $\langle t_i,t_{i+1}\rangle$.

\EXAMPLE {}

Let us suppose, that the
rate of unit of time inflation was
$0{.}1$
in time $0$ 
and rate of unit of time inflation was
$0{.}2$
in time $1$.
Then the rate of inflation should be changed in time 
$\langle 0,1\rangle$ between booth values 
for instance in this four following ways: 
$$
{\iota_{1}} :=\cases{ 0{.}1, &if $x<1$\cr
0{.}2, & if $x\geq 1$}
%$%
,\quad
%$
{\iota_{2}} := \frac {u^{2}}{10} + 0{.}1
%$%
,\quad
%$
{\iota_{3}}:= \frac {u}{10} + 0{.}1
%$%
,\quad
%$
{\iota_{4}} := \cases{ 0{.}1, &if $x\leq 0$\cr
0{.}2, & if $x> 0$}
\NEq$$%
The questions is:
What is the real value $X(1)$ %of  capital 
in time $1$  if the real value $X(0)$ of it 
in time $0$ was $100$?

In general, we have 
$$
X(1)={X(0)\over 
\hbox{e}^{\left(\int_0^1(\ln(1+\iota(u)\mathop{\rm d}u\right)}} 
\NEq$$

So, in our four cases we obtain:  
$$
\eqalign{
\iota=\iota_1=0{.}1\colon&\,
X(1)={\displaystyle \frac {100}{e^{
 \left(  \! \int _{0}^{1}{\ln}(1{.}1)\,\mathop{\rm d}u \!  \right) }}} =
90{.}9090909\dots
={100\over 1+0{.}1}
%$$
\cr
%$$
\iota=\iota_2={\displaystyle u\mapsto \frac {1}{10}} \,u^{2} + 0{.}1\colon &\,
X(1)={\displaystyle \frac {100}{e^{ \left(  \! \int _{0}^{1}{\ln}(1{.}1 + 
u^{2}/10)\,\mathop{\rm d}u \! \right) }}} =88{.}26551047\dots
%$$
\cr
%$$
\iota=\iota_3={u\mapsto\displaystyle \frac {1}{10}} \,u + 0{.}1\colon &\, 
X(1)={\displaystyle \frac {100}{e^{ \left(  \! \int _{0}^{1}{\rm ln}(
1{.}1 + u/10)\,\mathop{\rm d}u
 \!  \right) }}} =86{.}98393756\dots
%$$
\cr
%$$
\iota=\iota_4=0{.}2\colon &\, 
X(1)={\displaystyle \frac {100}{e^{
 \left(  \! \int _{0}^{1}{\rm ln}(1{.}2)\,\mathop{\rm d}u
 \!  \right) }}} =83{.}33333333\dots=
{100\over 1+0{.}2}
}
$$

\PROOF 

\noindent{\bf 1)}
The equation 
{%
$ \frac {{\rm x}(0)}{{\rm x}(t)}=e^{ \left(  \! 
\int _{0}^{t}{\rm ln}(1 + \iota (u))\,\mathop{\rm d}u
 \!  \right) }$%
}
gives the supposed results with the constant inflation:  
let us suppose, that 
$\iota (t)=I$
is constant:
$$
{\displaystyle \frac {{\rm X}(0)}{{\rm X}(t)}} 
=e^{ \left(  \! \int _{0}^{t}{\rm ln}(1 + I)\,\mathop{\rm d}u
 \!  \right) }
=e^{ \left(   {t}\cdot{\ln}(1 + I) \right) }
=e^{ \left({\ln}(1 + I)^t \right) }
=(1 + I)^t 
\NEq$$
q. e. d.

\noindent{\bf 2)}
Now let us suppose, that $\iota$ is piecewise constant and that it
has value
${I_{i}}$
 in every point of the
interval
${I_{i}}=( {t_{i}}, \,{t_{i + 1}}) \,;\,i=0 \dots n$.
Let
$\chi_A$
be the characteristic function of the set A, 
then $\iota(t)=\sum\chi_{(t_i,t_{i+1})}\cdot I_i$
$$
{\displaystyle \frac {{\rm X}(0)}{{\rm X}(t)}} 
=e^{ \left(  \! \int _{0}^{t}{\rm ln}(1 + 
\chi_{(t_i,t_{i+1})}\cdot I_i
\,\mathop{\rm d}u
 \!  \right) }
%$$
=
%$$
e^{ \left(  \! \sum _{i=0}^{n - 1}\,({t_{i + 1}}\,{\ln}(1 + {
I_{i}}) - {t_{i}}\,{\ln}(1 + {I_{i}})) \!  \right) }
\NEq$$
We shall use
$$
\displaylines{
e^{\left(\sum _{i=0}^{n - 1}\,{\zeta _{i}}\right)}=\prod _{i=
0}^{n - 1}\,e^{{\zeta _{i}}}
\cr
\leftline{hence}\cr
e^{ \left(  \! \sum _{i=0}^{n - 1}\,({t_{i + 1}}\,{\rm ln}(1 + {I
_{i}}) - {t_{i}}\,{\rm ln}(1 + {I_{i}})) \!  \right) }=
{\displaystyle \prod _{i=0}^{n - 1}} \,(1 + {I_{i}})^{({t_{i + 1}
} - {t_{i}})}
\cr
}
$$
consequently  
$$
{\displaystyle \frac {{\rm X}(0)}{{\rm X}(t)}} 
={\displaystyle \prod _{i=0}^{n - 1}} \,e^{({t_{i + 1}}
 - t_{i})
\,{\rm ln}
(1 + {I_{i})}}
=
{\displaystyle \prod _{i=0}^{n - 1}} \,e^{{\rm ln}
(1 + {I_{i})^{(t_{i + 1}
 - t_{i})}
}}
%$$
=
%$$
%{\displaystyle \frac {{\rm x}(0)}{{\rm x}(t)}} =
{\displaystyle \prod _{i=0}^{n - 1}} \,(1 + {I_{i}})^{({t_{i + 1
}} - {t_{i}})}
\NEq$$
q. e. d.

\K {Example with the real data of  index of prices of non-regulated prices in Czech 
republic}.

We choose the year as a unity of time. We  measured the index of non-regulated prices in 
every moment $P=\left\{1993+{i\over 12}\right\}$,
where $i$ is a natural number less or equal to $120$. The values  are relative. 
As a result of the measurement, \CPI is a function cutting the points:
\newcount\TMPnrA
\newcount\TMPnrB
\def\prepocet#1.#2#3{\TMPnrA=#1\relax\TMPnrB=#2#3\relax
\ifcase\TMPnrB 
\or \def\cast{+{1\over 12}}%
\or \def\cast{+{1\over 6}}%
\or \def\cast{+{1\over 4}}%
\or \def\cast{+{1\over 3}}%
\or \def\cast{+{5\over 12}}%
\or \def\cast{+{1\over 2}}%
\or \def\cast{+{7\over 12}}%
\or \def\cast{+{2\over 3}}%
\or \def\cast{+{3\over 4}}%
\or \def\cast{+{5\over 6}}%
\or \def\cast{+{11\over 12}}%
\or \advance \TMPnrA by 1 \def\cast{}\fi
$\number\TMPnrA \cast$}
\def\radek#1&#2\\{%
{\strut(\expandafter\prepocet#1;#2),~}}
{\let\par\radek
\advance\baselineskip by 2 pt
\obeylines
 1993.01& 91.76\\
 1993.02& 92.97\\
 1993.03& 93.54\\
 1993.04& 93.94\\
 1993.05& 94.32\\
 1993.06& 94.62\\
 1993.07& 95.42\\
 1993.08& 96.10\\
 1993.09& 97.53\\
 1993.10& 98.62\\
 1993.11& 99.17\\
 1993.12& 100.0\\
 1994.01& 100.76\\
 1994.02& 101.12\\
 1994.03& 101.40\\
 1994.04& 101.88\\
 1994.05& 102.27\\
 1994.06& 103.45\\
 1994.07& 103.81\\
 1994.08& 104.62\\ 
 1994.09& 106.13\\
 1994.10& 107.45\\
 1994.11& 108.46\\
 1994.12& 109.23\\
 1995.01& 110.57\\
 1995.02& 111.62\\
 1995.03& 111.97\\
 1995.04& 112.72\\
 1995.05& 113.26\\
 1995.06& 114.13\\
 1995.07& 113.47\\
 1995.08& 113.41\\
 1995.09& 114.41\\
 1995.10& 115.20\\
 1995.11& 116.09\\
 1995.12& 116.82\\
 1996.01& 118.97\\
 1996.02& 119.62\\
 1996.03& 120.43\\
 1996.04& 121.15\\
 1996.05& 121.96\\
 1996.06& 123.06\\
 1996.07& 123.18\\
 1996.08& 122.67\\
 1996.09& 123.04\\
 1996.10& 123.76\\
 1996.11& 124.44\\
 1996.12& 125.23\\
 1997.01& 126.28\\
 1997.02& 126.71\\
 1997.03& 126.85\\
 1997.04& 127.45\\
 1997.05& 127.60\\
 1997.06& 129.42\\
 1997.07& 129.67\\
 1997.08& 130.77\\
 1997.09& 131.58\\
 1997.10& 132.33\\
 1997.11& 133.01\\
 1997.12& 133.75\\
 1998.01& 135.76\\
 1998.02& 136.71\\
 1998.03& 136.85\\
 1998.04& 137.12\\
 1998.05& 137.26\\
 1998.06& 137.81\\
 1998.07& 137.53\\
 1998.08& 137.12\\
 1998.09& 137.25\\
 1998.10& 136.84\\
 1998.11& 136.43\\
 1998.12& 136.02\\
 1999.01& 136.70\\
 1999.02& 136.57\\
 1999.03& 136.29\\
 1999.04& 136.84\\
 1999.05& 136.70\\
 1999.06& 136.98\\
 1999.07& 136.98\\
 1999.08& 137.11\\
 1999.09& 136.98\\
 1999.10& 136.98\\
 1999.11& 137.39\\
 1999.12& 138.21\\
 2000.01& 139.04\\
 2000.02& 139.32\\
 2000.03& 139.32\\
 2000.04& 139.32\\
 2000.05& 139.74\\
 2000.06& 140.71\\
 2000.07& 141.42\\
 2000.08& 141.70\\
 2000.09& 141.56\\
 2000.10& 141.98\\
 2000.11& 142.13\\
 2000.12& 142.41\\
 2001.01& 143.26\\
 2001.02& 143.26\\
 2001.03& 143.26\\
 2001.04& 143.84\\
 2001.05& 144.99\\
 2001.06& 146.87\\
 2001.07& 147.90\\
 2001.08& 147.46\\
 2001.09& 145.98\\
 2001.10& 145.84\\
 2001.11& 145.69\\
 2001.12& 145.98\\
 2002.01& 147.30\\
 2002.02& 147.30\\
 2002.03& 147.0\\
 2002.04& 147.30\\
 2002.05& 147.15\\
 2002.06& 146.71\\
 2002.07& 147.15\\
 2002.08& 146.85\\
 2002.09& 145.82\\
 2002.10& 145.82\\
 2002.11& 145.68\\%
% 2002.12& 145.97\\
}
(2003; 147.97).
Now we are going to approximate the values of \CPI with the function. There are a lot of 
ways of
how to do it. We shall measure the accuracy  of approximation as the sum of 
squares of
differences between values of function and measured values in points where the measurement 
is made. Trivial approximation  is to put together  points by abscissas.
We obtain piecewise affine function:
and the accuracy  is equal to 0.
In this case the rate of the year inflation $\iota$  is the function:
{
\catcode`\[=\active
\def[#1,t<#2,#3<t]{#1&$t\in\langle#3,#2\rangle$\cr}
$$
\iota(t)= 
\cases{
[{e^{ 14.52058082\, \left(  14.52058082\,t- 28847.75831 \right) ^{-1}}
}-1,t< 1993.08333, 1993.00000<t]
[{e^{ 6.840273611\, \left(  6.840273611\,t- 13540.26521 \right) ^{-1}}
}-1,t< 1993.16666, 1993.08333<t]
[{e^{ 4.799616031\, \left(  4.799616031\,t- 9472.894168 \right) ^{-1}}
}-1,t< 1993.25000, 1993.16666<t]
[{e^{ 4.560182407\, \left(  4.560182407\,t- 8995.643826 \right) ^{-1}}
}-1,t< 1993.33333, 1993.25000<t]
[{e^{ 3.600144006\, \left(  3.600144006\,t- 7081.966879 \right) ^{-1}}
}-1,t< 1993.41666, 1993.33333<t]
[{e^{ 9.599232061\, \left(  9.599232061\,t- 19040.64915 \right) ^{-1}}
}-1,t< 1993.50000, 1993.41666<t]
[{e^{ 8.160326413\, \left(  8.160326413\,t- 16172.19129 \right) ^{-1}}
}-1,t< 1993.58333, 1993.50000<t]
[{e^{ 17.16068643\, \left(  17.16068643\,t- 34115.15901 \right) ^{-1}}
}-1,t< 1993.66666, 1993.58333<t]
[{e^{ 13.07895368\, \left(  13.07895368\,t- 25977.54380 \right) ^{-1}}
}-1,t< 1993.75000, 1993.66666<t]
[{e^{ 6.600264011\, \left(  6.600264011\,t- 13060.65643 \right) ^{-1}}
}-1,t< 1993.83333, 1993.75000<t]
[{e^{ 9.960398416\, \left(  9.960398416\,t- 19760.20401 \right) ^{-1}}
}-1,t< 1993.91666, 1993.83333<t]
[{e^{ 9.120364815\, \left(  9.120364815\,t- 18085.24781 \right) ^{-1}}
}-1,t< 1993.99999, 1993.91666<t]
[{e^{ 4.320172807\, \left(  4.320172807\,t- 8513.664947 \right) ^{-1}}
}-1,t< 1994.08332, 1993.99999<t]
[{e^{ 3.359731222\, \left(  3.359731222\,t- 6598.464123 \right) ^{-1}}
}-1,t< 1994.16666, 1994.08332<t]
[{e^{ 5.760230409\, \left(  5.760230409\,t- 11385.45902 \right) ^{-1}}
}-1,t< 1994.24999, 1994.16666<t]
[{e^{ 4.680187207\, \left(  4.680187207\,t- 9231.582863 \right) ^{-1}}
}-1,t< 1994.33332, 1994.24999<t]
\vdots
}
$$
}
and
the rate of inflation per the time interval
$\langle 1993{.}5,1994\rangle$
equals
$$
{\CPI(1994)\over\CPI(1993{.}5)}-1=
{e^{\int _{1993{.}5}^
{1994}\!\ln  \left(  1+{\iota} \left( z \right)  \right) \mathop{\rm d}z
}}
=0{.}05596
\NEq$$

We can use more sophisticated approximation, for instance:
let us approximate a trend by function
$
151.558469453+53.8746490595\cdot(5/12)^{(x-1992)}-108.487078769\cdot(3/4)^{(x-1992)}
$
which is the best approximation of the measured values by means of pair of functions out 
of the set $\left\{
(i/12)^{(x-1992)}\right.$;
$(i/12)^{(x-1993)}$;
$(x-1992)^{(i/12)}$;
$(x-1993)^{(i/12)}$;
$\ln(x-1992)$;
$\left.\ln(x-1991)
\right\}_{i=1}^{24}$
and constant
(accuracy (the sum of squares  of distances between values of the functions and the
measured values  equals  $237{.}324366060$))
Then let us approximate the rest by linear combination of functions
$$
\left\{
\sin(\pi   x/i);
\cos(\pi    x/i);
x\sin(\pi  x/i);
x\cos(\pi  x/i)\right\}_{i=1}^5
\NEq$$
Let  us omit 15 of the functions from this set of functions, the absence of which brings
the least loss of accuracy. We obtain the following approximation of CPI:
$
\CPI(t)=
151{.}558469453+53{.}8746490595\cdot(5/12)^{(x-1992)}-108{.}487078769\cdot(3/4)^{(x-1992)}
-0{.}08688205539+438{.}678532386\cdot \sin(1/5\cdot x\cdot \pi )-0{.}219746118617\cdot x\cdot \sin(1/5\cdot x\cdot \pi )
-0{.}392448110715\cdot x\cdot \cos(1/2\cdot x\cdot \pi )+0{.}000156347532197\cdot x\cdot 
\sin(x\cdot \pi )+782{.}647022069\cdot \cos(1/2\cdot x\cdot \pi ) 
$
(accuracy (the sum of squares  of distances between values of the functions and the 
measured values  equals  
$84{.}8709833772$))
Then the rate of inflation per year has in the time $x$ the value
$$\displaylines
{e^{%\left( 
53.8746490595\,{A\over B}}\cr
A=
\left ({\frac {5}{12}}\right )^{x-1992}\ln ({\frac {5}{12}})
- 108.487078769\,\left (3/4\right )^{x-1992}\ln (3/4)
+ 87.73570648\,\cos(1/5\,x\pi )\pi\cr 
- 0.219746118617\,\sin(1/5\,x\pi )
- 0.04394922372\,x\cos(1/5\,x\pi )\pi \cr
- 0.392448110715\,\cos(1/2\,x\pi )
+ 0.1962240554\,x\sin(1/2\,x\pi )\pi \cr
+ 0.000156347532197\,\sin(x\pi )
+ 0.000156347532197\,x\cos(x\pi )\pi \cr
- 391.3235111\,\sin(1/2\,x\pi )
\pi %\right)
\cr
B=
%\left ( 
151.4715874+ 53.8746490595\,\left ({\frac {5}{12}}\right )^{x-1992}
- 108.487078769\,\left (3/4\right )^{x-1992}
+ 438.678532386\,\sin(1/5\,x\pi )\cr
- 0.219746118617\,x\sin(1/5\,x\pi )
- 0.392448110715\,x\cos(1/2\,x\pi )\cr
+ 0.000156347532197\,x\sin(x\pi )
+ 782.647022069\,\cos(1/2\,x\pi )%\right )
^{-1}
}
$$
and
$$
{\CPI(1994)\over\CPI(1993{.}5)}-1=
e^{\int _{1993{.}5}^
{1994}\!\ln  \left(  1+{\iota} \left( z \right)  \right) \mathop{\rm d}z
}
=%0{.}05596
  0{.}04737
\NEq$$

But we can approximate the rate of inflation directly. If $\left(t_{n_i}\right)_{i\in J}$, 
$i<j\Rightarrow t_i<t_j$ are the moments, where  \CPI is measured, then we have to fulfil 
the equation: $$
\forall i\in J\colon \int_{t_i}^{t_{i+1}} \ln 
\left(1+\iota(t)\right)
\mathop{\rm d}t
=\ln\circ\CPI(t_{n+1}) -\ln\circ\CPI(t_{n})
\NEq
$$
\PripravCitaci{FunkcEq}
and the accuracy  of approximation should be measured as the measure of 
unreliableness of the equation.

More precisely:
Let us choose  the function:
$$\displaylines{
\iota\colon x\longmapsto
0.04944210+ 0.009094682\,\cos \left( x\pi  \right) + 
0.01215932\,\sin \left( 2\,x\pi  \right) +\hfill\cr\hfill+ 
0.01823419\,\sin \left( 4\,x\pi  \right) + 
0.05263253\,\sin \left( 1/3\,x\pi  \right) - \Nq\cr\hfill-
0.009667420\,\cos \left( 1/3\,x\pi  \right) - 
0.07212055\,\cos \left( 1/4\,x\pi  \right) + 
0.02834145\,\sin \left( 1/2\,x\pi  \right)
%- 0.602964297863- 0.0391382355972\,
%\cos \left( 1/2\,x\pi  \right) -\hfill\cr- 0.121546500735\,
%\cos \left( 1/3\,x\pi  \right) -0.263374831471\,
%\sin \left( 1/5\,x\pi  \right) -\cr\hfill\hfill- 0.634898132119\,
%\ln  \left( x-1992 \right) + 1.00628894644\, 
%\left( {\frac {13}{12}} \right) ^{x-1992}
.}$$
The question is, how accurately this function corresponds to the measured values. 
The interesting thing is, that the function do not approximate measured 
values,
but some other values, which depend  on the measured ones. And the attitude towards this 
two values is given by
(see. (\Cituj{FunkcEq})). 
So if we can compare the exactness of approximation with the same measure as we did in 
previous cases, we have to compute the number 
$$\sum_{i=1}^{120}\left(e^{\int_{t_i}^{t_{i+1}} \ln
\left(1+\iota(t)\right)\mathop{\rm d}t}-{\CPI(t_{n+1})\over\CPI(t_{n})}\right)^2\NEq$$ 
and this is the only relevant measure.
In our case we obtain:
 $$
4.9914
\NEq.$$
We used  the simplest function for the approximation of rate of 
inflation  we considered yet, and it gives us far the 
best approximation we have obtained!

\K {Conclusion}

Our sketch brings an open problem:  
what does the space of solution of (\Cituj{FunkcEq}) look like
and how can it help us to understand what kind of measure is this one which we called rate 
(of interest or inflation\dots).
It shows us the possibility approximate the measure of inflation except of index of prices 
and it suggests, that the solution should be better.

\let\ZacatekSkupiny{
\let\KonecSkupiny}

\edef\zobak{$>$}
{
\catcode`\>=\active
\gdef\Odrazka{%
\catcode`\>=\active
\let>\zobak}
}

\def\Maple{\ZacatekSkupiny\medskip
\parindent=0pt
\obeylines
\font\1=cmtt10\1
\catcode`\>=\active
\catcode`\^=11
\catcode`\_=11
\catcode`\#=11
\Odrazka
\Jadro}
\def\Jadro#1{#1\medskip\KonecSkupiny}

\K {Apendix}
Maple programs for computation.
Suppose that in iota we have function for approximation, in our case:
\Maple{
> iota:=.4944210e-1 + %
.9094682e-2 * cos(x*Pi) + .1215932e-1 * sin(2*x*Pi) + .1823419e-1 * sin(4*x*Pi)+ %
.5263253e-1 * sin(1/3*x*Pi) - .9667420e-2 * cos(1/3*x*Pi) - .7212055e-1 * cos(1/4*x*Pi)+ %
.2834145e-1 * sin(1/2*x*Pi);
}

Suppose, that values of \CPI are in the file Values and the corresponding moments of time, 
in which are values measured, are in file Points. Then we can compute the preciseness of 
approximation using following program:

\Maple{
> Mistake := 0; 
> for i to NumberOfPoints do 
> A := evalf(int(ln(1+iota),x = Points[1] .. Points[i+1])); 
> A := Values[1]*exp(A); 
> Mistake := Mistake+(A-Values[i+1])^2 
> end do:
> print(Sum('Delta'^2,i)=Mistake);
}

\newcount\Nr\Nr=1
\REFERENCES
\everypar{[\number\Nr] \advance\Nr by 1\space }

Jitka DupaŸov , Jan Hurt, Josef æœep n: Stochastic Modeling in economics
and
Finance, Kluwert Academic Publishers, 2002, ISBN 1-4020-0840-6.

Jean  Dieudonn': 
Treatise on analysis, {Vol. III},
Pure and Applied Mathematics, {Vol. 10-III}, Academic Press. New
York -- London, 1972, ISBN

J. Acz'l: Lectures on functional Equations and Their Applications, New York, Academic 
Press 1966

\end